\pdfoutput=1
\RequirePackage{ifpdf}
\ifpdf 
\documentclass[pdftex]{sigma}
\else
\documentclass{sigma}
\fi

\DeclareMathOperator{\Aut}{Aut}			
\newcommand{\A}{\mathcal{A}}

\newcommand{\Coo}{C^\infty}

\newcommand{\N}{\mathbb{N}}
\newcommand{\R}{\mathbb{R}}
\newcommand{\Z}{\mathbb{Z}}
\newcommand{\dd}{\mathrm{d}}
\newcommand{\NCA}{\A_{\theta}}

\begin{document}

\allowdisplaybreaks

\renewcommand{\thefootnote}{$\star$}

\renewcommand{\PaperNumber}{029}

\FirstPageHeading

\ShortArticleName{On Projections in the Noncommutative 2-Torus Algebra}

\ArticleName{On Projections in the Noncommutative\\
2-Torus Algebra\footnote{This paper is a~contribution to the Special Issue on Noncommutative Geometry and Quantum Groups
in honor of Marc A.~Rief\/fel.
The full collection is available at
\href{http://www.emis.de/journals/SIGMA/Rieffel.html}{http://www.emis.de/journals/SIGMA/Rieffel.html}}}

\Author{Micha{\l} ECKSTEIN}

\AuthorNameForHeading{M.~Eckstein}

\Address{Faculty of Mathematics and Computer Science, Jagellonian University,\\
ul.~{\L}ojasiewicza 6, 30-348 Krak\'ow, Poland}
\Email{\href{mailto:michal.eckstein@uj.edu.pl}{michal.eckstein@uj.edu.pl}}

\ArticleDates{Received December 09, 2013, in f\/inal form March 16, 2014; Published online March 23, 2014}

\Abstract{We investigate a~set of functional equations def\/ining a~projection in the noncommutative 2-torus algebra
$A_{\theta}$. The exact solutions of these provide various generalisations of the Powers--Rief\/fel projection.
By identifying the corresponding $K_0(A_{\theta})$ classes we get an insight into the structure of projections in
$A_{\theta}$.}

\Classification{46L80; 19A13; 19K14; 46L87}

\renewcommand{\thefootnote}{\arabic{footnote}}
\setcounter{footnote}{0}

\section{Introduction}

Projections (i.e.~selfadjoint, idempotent elements) in associative $*$-algebras are the main building blocks of the
algebraic $K$-theory.
Commutative $C^*$-algebras, which by Gelfand--Naimark theorem are equivalent to locally compact Hausdorf\/f spaces, do not
contain non-trivial projections, when the corresponding space is connected.
To determine the $K_0$ group of a~unital $C^*$-algebra~$\A$ one thus has to study the equivalence classes of projections
in the matrix algebra~$M_{\infty}(\A)$.
However, when one abandons the assumption of commutativity of the algebra one may encounter various non-trivial
projections in the algebra itself which, in some cases, are suf\/f\/icient to fully determine the group~$K_0(\A)$.

The $K$-theory of the noncommutative 2-torus algebra $A_{\theta}$, known also as the irrational rotation algebra, has
been thoroughly investigated in the 1980's.
From the works of Pimsner, Voiculescu and Rief\/fel (see~\cite{PimVoic,Rieffel2} and references therein) we know that
$K_0(A_{\theta}) \cong \Z \oplus \theta\Z \cong \Z \oplus \Z$.
In the case of noncommutative tori it turns out that projections in the algebra $A_{\theta}$ itself generate the whole
group $K_0(A_{\theta})$ (see \cite[Corollary~7.10]{Rieffel3}).
The $K_0$ class of a~projection is uniquely determined by its algebraic trace, so any two projections with the same
trace must be unitarily equivalent in $M_{\infty}(A_{\theta})$ (see \cite[Corollary 2.5]{Rieffel2}).
On the other hand, it has been already pointed out by Rief\/fel in~\cite{Rieffel3} that the structure of projections in
$A_{\theta}$ is more robust than it would appear from the $K$-theory level.

The purpose of this paper is to look closer into the structure of projections in $A_{\theta}$ itself.
Our main results are summarised in Theorems~\ref{lm1} and~\ref{lm2} in Section~\ref{general}.
Proposition~\ref{PropSym} addresses the problem of the existence of projections invariant under the f\/lip
automorphism~\eqref{flip}.
The statements are proven by an explicit construction of the relevant projections.
The latter may be useful in the applications where explicit formulae for projections are needed.
For an example see~\cite{exel}, where projections in $M_2(C(\mathbb{T}^2))$ were constructed using Rief\/fel's method,
which we generalise here\footnote{We thank the anonymous referee for pointing out this possible application to us.}.

The uses of noncommutative tori in physics are multifarious.
A most natural one concerns gauge theories developed in terms of f\/initely generated projective modules, which are
noncommutative counterparts of vector bundles~\cite{Connes,Varilly}.
Recently, the projections in the noncommutative torus algebra $A_{\theta}$ gained more interest in the context of string
theory~\cite{ConnesM,Witten}.
They turned out to be extrema of the tachyonic potential providing solitonic f\/ield solutions interpreted in terms of
D-branes~\cite{PRD,Krajewski,NFT}.
Moreover, the projections in $A_{\theta}$ are extensively used in the context of quantum anomalies~\cite{Magisterka,
Perrot}, knot theory~\cite{Knot} or theoretical engineering~\cite{Gabor}.

The paper is organised as follows: Below we recall some basic def\/initions to f\/ix notation and make the paper
self-contained.
In the next section we present a~set of functional equations def\/ining a~projection in $A_{\theta}$ and comment on the
adopted method of solving these.
We also compute the Chern class of a~projection satisfying these equations.
In Section~\ref{sol} we investigate some special solutions -- the Powers--Rief\/fel type projections.
These will serve us f\/irstly to provide an alternative proof of the important Corollary~7.10 from~\cite{Rieffel3}.
Secondly, we will use them as a~starting point for generalisations to come in Section~\ref{general}.
Section~\ref{symmetric} contains the discussion of the f\/lip-invariance~\eqref{flip}.
Finally, in Section~\ref{conc}, we make conclusions and discuss some open questions that arose during the presented
analysis.

The algebra of noncommutative 2-torus $A_{\theta}$ is the universal $C^*$-algebra generated by two unitaries $U$, $V$
satisfying the following commutation relation
\begin{gather}
\label{CommRel}
VU = e^{2 \pi i \theta} UV
\end{gather}
for some real parameter $\theta \in [0,1[$, which we assume to be irrational.

We shall work with $\NCA$, a~dense $*$-subalgebra~\cite{Masoud,Varilly} of $A_{\theta}$ which is made up of ``smooth''
elements of the form
\begin{gather}
\label{a}
A_{\theta} \supset \NCA \ni a= \sum\limits_{(m,n)\in\Z^2} a_{m,n} U^m V^n,
\qquad
\text{with}
\quad
\{a_{m,n}\} \in \mathcal{S}\big(\Z^2\big),
\end{gather}
where $\mathcal{S}(\Z^2)$ denotes the space of Schwartz sequences on $\Z^2$, i.e.~
\begin{gather*}
\{a_{m,n}\} \in \mathcal{S}\big(\Z^2\big)
\quad
\Longleftrightarrow
\quad
\sup_{m,n \in \Z} \big(1+ m^2 + n^2\big)^k \vert a_{m,n}\vert < \infty,
\qquad
\text{for all}
\quad
k \in \N.
\end{gather*}
Let us note that $\NCA$ is a~Fr\'echet pre-$C^*$-algebra~\cite{Varilly}, hence $K_0(\NCA) \cong K_0(A_{\theta})$ as
Abelian groups~\cite[Theorem 3.44]{Elements}.

For the purposes of this paper we will only consider elements of $\NCA$ with $a_{m,n} \neq 0$ for a~f\/inite range of the
index $m$.
These can be written as
\begin{gather}
\label{b}
b = \sum\limits_{m\text{-f\/inite}} \sum\limits_{n \in \Z} b_{m,n} U^m V^n,
\qquad
\text{with}
\quad
\{b_{m,n}\}_{n \in \Z} \in \mathcal{S}(\Z)
\quad
\text{for any}~m.
\end{gather}
Recall that any function $f \in C^{\infty}(S^1)$, regarded as a~periodic function on $\R$ with period 1, has a~Fourier
series presentation
\begin{gather*}
f(x) = \sum\limits_{n \in \Z} f_n e^{2\pi i x n},
\qquad
\text{with}
\quad
\{f_n\} \in \mathcal{S}(\Z).
\end{gather*}
Thus, via the funcional calculus, we may uniquely def\/ine an element from $\NCA$ by
\begin{gather*}
f(V) = \sum\limits_{n \in \Z} f_n V^n,
\qquad
\text{for}
\quad
f \in C^{\infty}(S^1)
\end{gather*}
and the expression~\eqref{b} is conveniently rewritten as
\begin{gather}
\label{b1}
b = \sum\limits_{m\text{-f\/inite}} U^m b_m(V),
\qquad
\text{with}
\quad
b_m \in C^{\infty}(S^1)
\quad
\text{for any}~m.
\end{gather}

The noncommutative 2-torus algebra is equipped with a~canonical trace~\cite{Masoud,Varilly}, which on the elements of
the form~\eqref{b1} can be expressed as (compare with~\cite[p.~415]{Rieffel1})
\begin{gather}
\label{Atr}
\tau(b) = b_{0,0} = \int_0^1 b_0(x) \, \dd x.
\end{gather}
We shall use it to determine the $K_0$ class of a~projection on the strength of Corollary~2.5 in~\cite{Rieffel2}.

To make a~connection with the original framework of Rief\/fel's construction~\cite{Rieffel1} we shall describe~$A_{\theta}$ as the crossed product algebra~$C(S^1) \rtimes \Z$.
Here $\Z$ acts on $S^1$ by rotations by $2 \pi \theta$ and thus induces an action of~$\Z$ as automorphisms of~$C(S^1)$.
A convenient concrete realisation of~$A_{\theta}$ as bounded operators on $L^2(S^1)$ is obtained with
\begin{gather*}
(U f)(x) = f(x+ \theta), \qquad (V f)(x) = e^{2\pi i x} f(x).
\end{gather*}
Note that the elements of the form~\eqref{b1} fall in the dense $*$-subalgebra of $A_{\theta}$ considered by Rief\/fel
in~\cite{Rieffel1}.
The only dif\/ference is that we have chosen to work with the smooth functions on $S^1$ rather than continuous ones.
This is necessary as we want to compute the Chern number~\cite{ConnesChern} of the relevant projections, which requires~$b_n$ functions to be dif\/ferentiable.

Let $\delta_1$, $\delta_2$ be the basic unbounded derivations of $A_{\theta}$, which act on the generators as
\begin{gather*}
\delta_1 U = 2 \pi i U, \qquad \delta_1 V = 0, \qquad \delta_2 U = 0, \qquad \delta_2 V = 2 \pi i V.
\end{gather*}
Then the Chern number of a~projection reads
\begin{gather*}
c_1(p) = \frac{1}{2 \pi i} \tau \big( p ( \delta_1 p \delta_2 p - \delta_2 p \delta_1 p ) \big).
\end{gather*}
The Chern number is related to the index of a~Fredholm operator and thus it is always an integer
(see~\cite[Theorem~11]{ConnesChern}).

\section[Equations for a projection in $\A_{\theta}$]{Equations for a projection in $\boldsymbol{\A_{\theta}}$}
\label{equations}

Having recalled the basic features of the noncommutative 2-torus algebra we are ready to investigate the structure of
projections in it.

Let us consider the following element of $\A_{\theta}$:
\begin{gather}
\label{GenP}
p = \sum\limits_{n = -M}^M U^n p_n(V), \qquad\text{for some} \quad M \in \N.
\end{gather}
The conditions for $p$ to be a~projection yield a~set of functional equations for the functions $p_i \in \Coo(\R/\Z)$
\begin{gather}
p_k(x) = \overline{p_{-k}(x+k\theta)}, \qquad\text{for}\quad k = -M, \ldots, M,
\label{selfad}
\\
p_k(x) = \sum\limits_{m,a = -M}^M p_m(x+a \theta)p_a(x)   \delta_{m+a,k}, \qquad\text{for}\quad k = -M, \ldots, M,
\label{idemp1}
\\
0 = \sum\limits_{m,a = -M}^M p_m(x+a \theta)p_a(x)  \delta_{m+a,k}, \qquad\text{for}\quad k < -M
\quad
\text{and}
\quad
k > M.
\label{idemp2}
\end{gather}
Some of the above equations are redundant and the number of independent ones is $3M +2$.
It can be easily seen by noticing that equations~\eqref{selfad}--\eqref{idemp2} with $k < 0$ are equivalent to those with
$k >0$, because the functions $p_k$ with negative indices are actually def\/ined by~\eqref{selfad} with $k > 0$.
For $M = 0$ formulae~\eqref{selfad}--\eqref{idemp2} imply $p_0(x) \equiv 1$ as one may expect.
When $M = 1$ one obtains the familiar Powers--Rief\/fel equations~\cite{Rieffel1}.
However, for $M \geq 2$ the equations become more and more involved and even the existence of a~solution is not obvious.
In~\cite{Magisterka} we found four particular solutions to~\eqref{selfad}--\eqref{idemp2} with $M = 2$, which represent
dif\/ferent classes of $K_0(\A_{\theta})$.
In the next sections we present a~generalisation of the construction given in~\cite{Magisterka}.
Before we start solving the equations~\eqref{selfad}--\eqref{idemp2} let us adopt the following def\/inition.

\begin{definition}
We say a~projection in $\NCA$ is of \emph{order} $M$ if it is of the form~\eqref{GenP} and $p_M \neq 0$.
\end{definition}

We shall not attempt to provide a~general solution to~\eqref{selfad}--\eqref{idemp2}, but rather present a~class of special
solutions.
Nevertheless, this class turns out to be large enough to accommodate the known projections as well as a~number of new
ones.

We will consider only real-valued functions although~\eqref{selfad} requires only $p_0$ to be real.
Moreover, we have already noted that~\eqref{selfad} def\/ines the functions $p_k$ for $k<0$ and it is convenient to get
rid of the functions $p_k$ with negative index $k$ in the equations~\eqref{idemp1} and~\eqref{idemp2} before solving
them.
Our special solutions will be such that each summand on the r.h.s.\ of~\eqref{idemp2} is equal to zero independently.
The same should hold for summands of~\eqref{idemp1} with $k>0$ excluding those with $m=0$ or $a=0$, these are combined
to form equations
\begin{gather}
\label{eqK}
p_k(x) \big( p_0(x) + p_0(x+k \theta) - 1 \big) = 0, \qquad \text{for}\quad k = 1, \ldots, M,
\end{gather}
which we also require to be satisf\/ied independently.
The equations~\eqref{idemp1} with $k<0$ are redundant and the case $k=0$ cannot be split into independent equations.
After~\eqref{selfad} is substituted into~\eqref{idemp1} for $k=0$ we obtain
\begin{gather}
p_M^2(x-M \theta) + p_M^2(x) + p_{M-1}^2(x-(M-1) \theta) + p_{M-1}^2(x)
\nonumber
\\
\qquad{}
+ \cdots + p_1^2(x-\theta) + p_1^2(x) + p_0(x) \big( p_0(x) - 1 \big) = 0.\label{eq0}
\end{gather}
In the forthcoming sections we provide a~systematic method of constructing projections of a~given trace~\eqref{Atr} and
order, that will satisfy the equations~\eqref{selfad}--\eqref{idemp2} ref\/ined according to the above-listed conditions.

Before we start, let us compute the Chern number of a~projection of order $M$, as this quantity might prove useful in
the task of classif\/ication.
\begin{proposition}
The Chern number of a~projection $p$ of order $M$ reads
\begin{gather}
c_1(p) = \sum\limits_{n = 1}^{M} \sum\limits_{k = -M}^{M-n} \int_{0}^{1} {\rm d}x
\Big( n\, \overline{p_{k+n}(x)} \big[p_n(x + k \theta) p_k'(x) - p_n(x) p_k'(x+n \theta) \big]
\nonumber
\\
\phantom{c_1(p)=}
{}+ \big( \{k,n\} \longleftrightarrow \{-k,-n\} \big) \Big).
\label{ChernGen}
\end{gather}
Moreover, for any projection constructed with the method adopted in this paper, the formula~\eqref{ChernGen} simplifies to
\begin{gather}
\label{ChernSimp}
c_1(p) = 6 \int_0^1 {\rm d}x \sum\limits_{n=1}^{M} n \, p_n(x)^2 p_0'(x).
\end{gather}
\end{proposition}

\begin{proof}
The application of the derivations $\delta_1$, $\delta_2$ to a~projection $p$ of the form~\eqref{GenP} yields
\begin{gather*}
\delta_1 p = 2 \pi i \sum\limits_{n = -M}^M n U^n p_n(V), \qquad \delta_2 p = \sum\limits_{n = -M}^M U^n p_n'(V).
\end{gather*}
Let us denote $p_k:= 0$ for $\vert k \vert > M$.
Then using~\eqref{CommRel},~\eqref{Atr} and~\eqref{selfad} we obtain
\begin{gather*}
c_1(p) = \sum\limits_{j,k,n = -M}^M \tau \left( U^j p_j(V) \big[ k U^k p_k(V) U^n p_n'(V) - n U^k p_k'(V) U^n p_n(V)
\big] \right)
\\
\phantom{c_1(p)}
= \sum\limits_{j,k,n = -M}^M \tau \left( U^{j+k+n} p_j(e^{2\pi i \theta (k+n)} V) \big[ k p_k(e^{2\pi i \theta n} V)
p_n'(V) - n p_k'(e^{2\pi i \theta n}V) p_n(V) \big] \right)
\\
\phantom{c_1(p)}
= \sum\limits_{k,n = -M}^M \int_0^1 {\rm d}x \, p_{-k-n}(x+ (k+n)\theta) \big[ k p_k(x + \theta n) p_n'(x) - n p_k'(x + n
\theta) p_n(x) \big]
\\
\phantom{c_1(p)}
= \sum\limits_{k,n = -M}^M n \int_0^1 {\rm d}x \, \overline{p_{k+n}(x)} \big[ p_n(x + \theta k) p_k'(x) - p_k'(x + n
\theta) p_n(x) \big].
\end{gather*}
In the last equality we have relabelled the indices $k \leftrightarrow n$ in the f\/irst term in the square bracket.
The formula~\eqref{ChernGen} is just the above expression, with the fact $p_k = 0$ for $\vert k \vert > M$ taken into
account.

Now, in the method adopted in this paper we have assumed that the summands of the r.h.s.\ of~\eqref{idemp1}
and~\eqref{idemp2} are equal to zero independently, excluding those with $m=0$ or $a=0$.
This means that for the considered projections only the terms $k=0$ and $k=-n$ will contribute to the sum
in~\eqref{ChernGen} and we have
\begin{gather*}
c_1(p) = \sum\limits_{n=1}^{M} \int_0^1 {\rm d}x \, n \Big( p_n(x)^2 p_0'(x) - p_n^2(x) p_0'(x+n\theta) + 2 p_0(x) p_n(x-n
\theta) p_{-n}'(x)
\\
\phantom{c_1(p)=}
{}- 2 p_0(x) p_n(x) p_{-n}'(x+n \theta) + p_n(x-n\theta)^2 p_0'(x-n\theta) - p_n^2(x-n\theta) p_0'(x) \Big)
\\
\phantom{c_1(p)}
= 2 \sum\limits_{n=1}^{M} \int_0^1 {\rm d}x \, n \Big( p_n(x)^2 p_0'(x) - p_n^2(x) p_0'(x+n\theta)
\\
\phantom{c_1(p)=}
{}+ p_0(x) p_n(x-n \theta) p_{n}'(x-n\theta) - p_0(x) p_n(x) p_{n}'(x) \Big),
\end{gather*}
We have used the formula~\eqref{selfad} together with the assumption of all $p_n$ being real and the fact that shifting
the integration variable does not change the value of the integral.
Now, let us make use of the periodicity of the integrant to integrate by parts one of the terms and shift the variable
in another:
\begin{gather*}
c_1(p) = 2 \sum\limits_{n=1}^{M} \int_0^1 {\rm d}x \, n \Big( p_n(x)^2 p_0'(x) + 2 p_n(x) p_n'(x) p_0(x+n\theta)
\\
\phantom{c_1(p)=}
{}+ p_0(x+n \theta) p_n(x) p_{n}'(x) - p_0(x) p_n(x) p_{n}'(x) \Big).
\end{gather*}
Let us note that equations~\eqref{eqK} imply
\begin{gather*}
\int_0^1 {\rm d}x \, p_n'(x) p_n(x) p_0 (x+n\theta) = \int_0^1 {\rm d}x \, p_n'(x) p_n(x) - \int_0^1 {\rm d}x \, p_n'(x) p_n(x) p_0 (x)
\\
\hphantom{\int_0^1 {\rm d}x \, p_n'(x) p_n(x) p_0 (x+n\theta)}{}
= - \int_0^1 {\rm d}x \, p_n'(x) p_n(x) p_0 (x),
\end{gather*}
since $p_n'(x) p_n(x)$ is a~total derivative.
Finally, we obtain
\begin{gather*}
c_1(p) = 2\sum\limits_{n=1}^{M} \int_0^1 {\rm d}x \, n \left( p_n(x)^2 p_0'(x) - 4 p_0(x) p_n(x) p_{n}'(x) \right) = 6
\sum\limits_{n=1}^{M} \int_0^1 {\rm d}x \, n\, p_n(x)^2 p_0'(x).
\tag*{\qed}
\end{gather*}
\renewcommand{\qed}{}
\end{proof}

\section{Powers--Rief\/fel type projections}
\label{sol}

We start with recalling the construction of the Powers--Rief\/fel projection in a~slightly more general framework.
It will serve us as a~starting point for generalisations to come in the next section.

\begin{figure}[t]\centering
\includegraphics{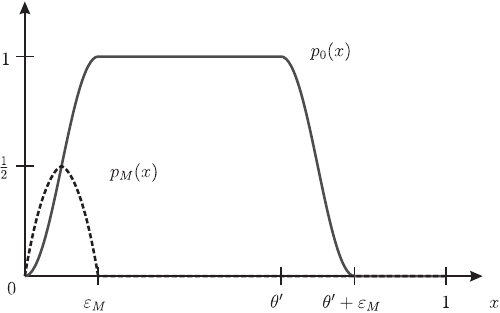}
\caption{Depiction of functions constituing a~Powers--Rief\/fel type projection. We have $ 0< \varepsilon_M \leq \theta' $, $\varepsilon_M + \theta' \leq 1$.}\label{f1}
\end{figure}

If one sets $p_k = 0$ for all $1 \leq k \leq M-1$ then~\eqref{selfad}--\eqref{idemp2} reduce to the Powers--Rief\/fel
equations with parameter $M \theta$ \cite[Theorem~1.1]{Rieffel1}
\begin{gather}
p_M(x+M \theta)p_M(x) = 0,
\label{PRM1}
\\
p_M^2(x) + p_M^2(x-M\theta) + p_0(x)(p_0(x) - 1) = 0,
\label{PRM2}
\\
p_M(x) \big( 1 - p_0(x) - p_0(x+ M\theta) \big) = 0.
\label{PRM3}
\end{gather}
A standard solution to~\eqref{PRM1}--\eqref{PRM3} is known as a~Powers--Rief\/fel type projection~\cite{EEproj,NFT}
\begin{gather}
p_0(x)=
\begin{cases}
d_M(x), & 0 \leq x < \varepsilon_M,
\\
1, & \varepsilon_M \leq x \leq M \theta,
\\
1 - d_{M}(x -M \theta), & M \theta \leq x < M \theta + \varepsilon_M,
\\
0, & M \theta + \varepsilon_M \leq x \leq 1,
\end{cases}
\label{PRT0}
\\
p_M(x)=
\begin{cases}
\sqrt{d_M(x)(1-d_M(x))}, & 0 \leq x < \varepsilon_M,
\\
0, & \varepsilon_M < x \leq 1,
\end{cases}
\label{PRTM}
\end{gather}
where $\theta' = M \theta - \lfloor M \theta \rfloor$ and $d_M$ is a~smooth function with $d_M(0) = 0$,
$d_M(\varepsilon_M) = 1$.
The functions $p_0$ and $p_1$ are depicted in Fig.~\ref{f1}.

Let us stress that we do not assume that the $d_M$ function starts growing directly at $x=0$ as shown on Fig.~\ref{f1}.
We may take $d_M$ such that $d_M = 0$ for $x \in [0,\delta_M]$ with some $\delta_M < \varepsilon_M$ and then smoothly
growing to reach 1 at $x=\varepsilon_M$.
This ensures that what we call here a~Powers--Rief\/fel type projection is suf\/f\/iciently general to incorporate the existing
def\/initions (see~\cite{NFT} for instance).

Let us now discuss the properties of these projections.
First of all, note that due to the periodicity of $p_i$ functions, equations~\eqref{PRM1}--\eqref{PRM3} are invariant with
respect to the transformation $M \theta \to M \theta + z$ for any $z \in \Z$.
This means that a~Powers--Rief\/fel type projection of order $M$ has the algebraic trace~\eqref{Atr} equal to $\theta' = M
\theta - \lfloor M \theta \rfloor$.
Since~$\theta$ is irrational we have inf\/initely many~$M$ such that
\begin{gather*}
0 < M \theta - n < 1 \quad
\Longleftrightarrow
\quad \frac{n}{M} < \theta < \frac{n+1}{M}.
\end{gather*}

Hence, the following proposition (which is also a~consequence of the Corollary~7.10 in~\cite{Rieffel3}) holds.
\begin{proposition}
\label{Thm}
The algebra $\A_{\theta}$ contains projections representing infinitely many different classes of $K_0(\A_{\theta})$.
\end{proposition}

Another point of view one may adopt for the projection~\eqref{PRT0}--\eqref{PRTM} is that for any f\/ixed $M$ it is the
standard Powers--Rief\/fel projection~\cite{Rieffel1} in the subalgebra of $\NCA$ generated by $U^M$ and $V$.
This fact may be used to construct an approximation of $\NCA$ in terms of two algebras of matrix valued functions on
$S^1$~\cite{EEproj,NFT}.

For $p^{[M]}$  a~Powers--Rief\/fel type projection of order $M$ the formula~\eqref{ChernSimp} gives
\begin{gather*}
c_1(p^{[M]}) = 6 M \int_{0}^{\varepsilon_M} {\rm d}x \, d_M(x) (1-d_M(x)) d_M'(x) = 6 M \left( \frac{d_M(x)^2}{2} -
\frac{d_M(x)^3}{3} \right) \Big\vert^{\varepsilon_M}_{0} = M.
\end{gather*}
This is in accordance with the result of~\cite{ConnesChern} stating that if $\tau(p) = \vert a-b \theta \vert$ then
$c_1(p) = \pm b$.

From the $K$-theoretic point of view, these projections are suf\/f\/icient to understand the structure of the equivalence
classes of projective modules over~$\A_{\theta}$.
On the other hand, the algebra~$\A_{\theta}$ contains other interesting projections, which we shall present in next
section.

\section[More general projections in $\NCA$]{More general projections in $\boldsymbol{\NCA}$}
\label{general}

Let us now see what kind of projections one can get by letting functions $p_k$ in~\eqref{GenP} to be non-zero for some
of the indices $k \in \{1, \ldots, M-1 \}$.
The results are summarised in the following Theorems.

\begin{theorem}\label{lm1}
A projection of order $M$ may represent the $K_0(\A_{\theta})$ class $[n \theta]$, as well as the class $[1 - n
\theta]$, for all $n = 1, 2, \ldots, \frac{1}{2}M(M+1)$, provided that $0< \theta < 1/\max(n,M)$.
\end{theorem}

By $[n\theta] \in K_0(\NCA)$ we denote the $K_0$ class represented by a~projection $p \in \NCA$ with $\tau(p) = n
\theta$.

\begin{theorem}\label{lm2}
The equations~\eqref{selfad}--\eqref{idemp2} for a~projection of order $M$ admit solutions with $p_k \neq 0$ for every $k
\in \{0, \ldots, M\}$ whenever $0< \theta < 1/M$.
\end{theorem}

We shall start with the proof of Theorem~\ref{lm1} by showing how to use the functions $p_k$ to increase or decrease the
trace of a~Powers--Rief\/fel type projection.
Then we present a~method of including the remaining $p_k$ functions into the projections constructed in the previous
proof without changing its traces.
In this way we will prove Theorem~\ref{lm2}.
Both proofs are constructive so we are able to plot some examples of the $p_0$ functions of the relevant projections
which, as we shall see, determine all of the other functions $p_k$ for $k \neq 0$.
A brief discussion of the assumptions limiting the $\theta$ parameter may be found in Section~\ref{conc}.

\begin{proof}
[Proof of Theorem~\ref{lm1}] Let us start with the case of $\tau(p) = n \theta > M \theta$.
We shall begin with a~Powers--Rief\/fel type projection as def\/ined in~\eqref{PRT0}--\eqref{PRTM}.
First note that if $M \theta < 1$ then the functions $p_0$ and $p_M$ of the Powers--Rief\/fel type projection of order $M$
vanish for $x \geq M \theta + \varepsilon_M$.
If $\theta$ is small enough (i.e.~$(M+k) \theta < 1$) then we can ``glue'' a~Powers--Rief\/fel type projection of trace $k
\theta$ to the previous one.
Namely, let us keep the def\/inition of $p_0$ on $[0,M \theta + \varepsilon_M]$ (see~\eqref{PRT0}) and set
\begin{gather*}
p_0(x) =
\begin{cases}
d_k(x), & M \theta + \varepsilon_M \leq x < M \theta + \varepsilon_M + \varepsilon_k,
\\
1, & M \theta + \varepsilon_M + \varepsilon_k \leq x \leq (M+k) \theta + \varepsilon_M,
\\
1 - d_{k}(x - k \theta), & (M+k) \theta + \varepsilon_M \leq x < (M+k) \theta + \varepsilon_M + \varepsilon_k,
\\
0, & (M+k) \theta + \varepsilon_M + \varepsilon_k \leq x \leq 1,
\end{cases}
\\
p_k(x) =
\begin{cases}
\sqrt{d_k(x)(1-d_k(x))}, & M \theta + \varepsilon_M \leq x < M \theta + \varepsilon_M + \varepsilon_k,
\\
0, & \text{elsewhere}
\end{cases}
\end{gather*}
for a~smooth function $d_{k}$ with $d_k(M \theta + \varepsilon_M) = 0$, $d_k(M \theta + \varepsilon_M + \varepsilon_k) =
1$ and a~small parame\-ter~$\varepsilon_k$.
The summands of~\eqref{idemp1} and~\eqref{idemp2}, which we have assumed to be equal to zero independently, have the
form $p_m(x+a\theta) p_a(x)$.
This means that all of the non-zero functions $p_k$ for $k \neq 0$ shifted to the interval $x \in [0,\theta]$ must not
intersect.
The latter can be fulf\/illed by restricting the parame\-ters~$\varepsilon$ such that
\begin{gather}
\label{eps2}
0< \varepsilon_M \leq M \theta, \qquad 0 < \varepsilon_k \leq k \theta, \qquad \varepsilon_M + \varepsilon_k + (M+k)\theta \leq 1
\end{gather}
implying that equation~\eqref{eq0} reduces to two equations of the form~\eqref{PRM2}.
Namely for $x \in [0,\varepsilon_M] \cup [M\theta,M\theta+\varepsilon_M]$ and for $ x \in
[M\theta+\varepsilon_M,M\theta+\varepsilon_M+\varepsilon_k] \cup
[(M+k)\theta+\varepsilon_M,(M+k)\theta+\varepsilon_M+\varepsilon_k]$ we have respectively
\begin{gather*}
p_0(x)(1-p_0(x)) = p_M^2(x) + p_M^2(x-M\theta),
\\
p_0(x)(1-p_0(x)) = p_k^2(x) + p_k^2(x-k\theta).
\end{gather*}
These equations are satisf\/ied by the construction of $p_k$ and $p_M$.
On the remaining part of the interval $[0,1]$ the equation~\eqref{eq0} is trivially satisf\/ied, since both l.h.s.\ and r.h.s.\
are equal to~0.
By the same argument, equation~\eqref{eqK} remains satisf\/ied, as it is satisf\/ied for both Powers--Rief\/fel type
projections independently.
Thus, we have obtained a~new projection with a~trace $(M+k)\theta$.
Examples of $p_0$ functions def\/ining such projections are depicted in Fig.~\ref{f2}.

If the parameter $\theta$ is small enough (i.e.~$n \theta < 1$) we can continue the process of ``glueing''
Powers--Rief\/fel type projections to obtain a~projection of trace $n \theta$, with $n \geq M$.
If one makes use of all of the functions $p_k$ with $1\leq k \leq M-1$ to increase the trace, one will end with
a~projection bearing the trace $(1 + 2 + \cdots + M) \theta = \frac{1}{2}M(M+1) \theta$.
The only thing one has to take care of are the conditions satisf\/ied by the parameters $\varepsilon_k$.
The restrictions~\eqref{eps2} may be easily generalised to the case of non-vanishing $p_{k_s}$ functions with $s\in[1,M-1]$:
\begin{gather}
0< \varepsilon_{k_j} \leq k_j \theta, \qquad\text{for}\quad 1 \leq j \leq s,
\nonumber
\\
\varepsilon_{k_1} + \cdots + \varepsilon_{k_s} + \varepsilon_M + n\theta \leq 1,
\qquad
\text{with}
\quad
n = k_1 + \cdots + k_s + M.
\label{eps}
\end{gather}

\begin{figure}[t]
\centering
\includegraphics{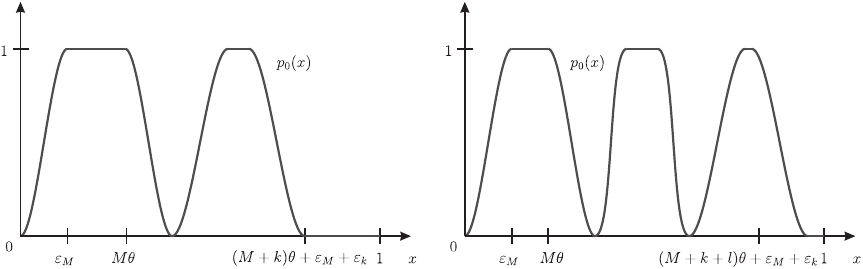}
\caption{Examples of $p_0$ functions for projections with traces $(M + k)\theta$ and $(M + k + l)\theta$.}  \label{f2}
\end{figure}

Let us note, that the above construction can be obtained (for $n\theta < 1$) by taking a~sum of~$s$ mutually orthogonal
Powers--Rief\/fel type projections $p^{[k_j]}$ of respective orders~$k_{j}$.
Indeed, one can easily check that the functional equations resulting form the projection and orthogonality conditions
\begin{gather*}
(p^{[k_i]})^2 = (p^{[k_i]})^* = p^{[k_i]}, \qquad p^{[k_i]} p^{[k_j]} = p^{[k_j]} p^{[k_i]} = 0,
\qquad\text{for}\quad
1 \leq i \neq j \leq s,
\end{gather*}
coincide with the ones derived in Section~\ref{equations}.
This process of ``glueing'' mutually orthogonal Powers--Rief\/fel type projections appeared already in~\cite{EEproj} and
was extensively used therein.
It has also been presented in~\cite{PRD} in a~more similar form to the one shown above.

Let us now consider the case of projections of order $M$ and trace $n \theta$ with $1 \leq n < M$.
Again, we shall use as a~starting point a~Powers--Rief\/fel type projection~\eqref{PRT0}--\eqref{PRTM}, but now we will
``cut out'' a~part of it.
Let us set
\begin{gather}
p_0(x)=
\begin{cases}
d_k(x), & \varepsilon_M \leq x < \varepsilon_M + \varepsilon_k,
\\
0, & \varepsilon_M + \varepsilon_k \leq x \leq k \theta + \varepsilon_M,
\\
1 - d_{k}(x - k \theta), & k \theta + \varepsilon_M \leq x < k \theta + \varepsilon_M + \varepsilon_k,
\\
1, & k \theta + \varepsilon_M + \varepsilon_k \leq x \leq M \theta,
\end{cases}
\label{cut0}
\\
p_k(x) =
\begin{cases}
\sqrt{d_k(x)(1-d_k(x))}, & \varepsilon_M \leq x < \varepsilon_M + \varepsilon_k,
\\
0, & \text{elsewhere}
\end{cases}
\label{cutk}
\end{gather}
with a~smooth function $d_k$ such that $d_k(\varepsilon_M) = 1$, $d_k(\varepsilon_M + \varepsilon_k) = 0$.
The conditions $0 < \varepsilon_M \leq M \theta $, $0 < \varepsilon_k \leq k \theta$ and $\varepsilon_M + M\theta \leq
1$ should be satisf\/ied.
The situation is now completely analogous to the case of ``glued'' projections and the same arguments apply.
A projection obtained in this way bears the trace $(M-k) \theta$ for $1 \leq k \leq M-1$ (see Fig.~\ref{f3}).

\begin{figure}[t]\centering
\includegraphics{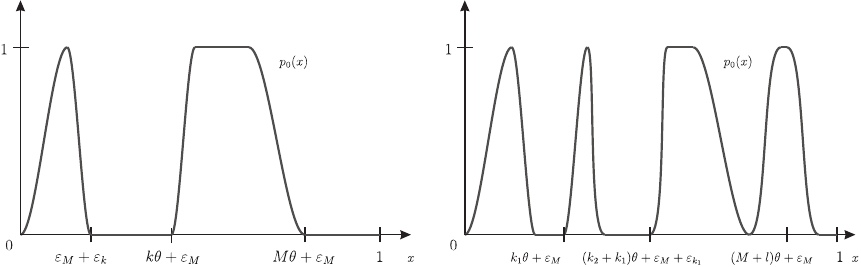}
\caption{Examples of $p_0$ functions for projections of trace $(M-k)\theta$ and $(M-k_1-k_2+l)\theta$.}\label{f3}
\end{figure}

To end the proof of Theorem~\ref{lm1} it remains just to recall that if $p$ is a~projection then obviously $1-p$ is so.
This means that all of the considerations hold for projections of traces $(1 - n \theta)$, one simply should take $1 -
p_0$ instead of $p_0$ and leave $p_k$ for $k \neq 0$ as they are.
\end{proof}

The presented proof provides a~great variety of possible projections with a~given trace, which have, in general,
dif\/ferent orders.
Let us notice that the two procedures of increasing and decreasing the trace of a~projection of a~given order can be
applied simultaneously and in arbitrary sequence (see Figs.~\ref{f3} and~\ref{f5}).
One only has to choose well the parameters $\varepsilon_{k}$ to have the equations~\eqref{eps} satisf\/ied.
These equations guarantee that the functions $d_{k}$ do not superpose and the
equations~\eqref{selfad}--\eqref{idemp2} remain satisf\/ied.
This leads to an enormous number of projections if the order $M$ is big enough.
Let us now pass on to the most general projections we were able to construct with the adopted method.

\begin{proof}[Proof of Theorem~\ref{lm2}]
In fact one can let all $p_k$ functions to be non-zero by incorporating to~$p_0$ some
``bump functions'' $d_k$.
As a~starting point, one should take an arbitrary projection def\/ined in Section~\ref{sol} or~\ref{general}.
For sake of simplicity let us now denote by $k$ a~free index, i.e.~we have $p_k = 0$ in our starting point projection.
Now, if one sets $p_0(x)=d_k(x)$ for $x\in[\delta_k,\delta_k+\varepsilon_k]$, with
$d_k(\delta_k)=d_k(\delta_k+\varepsilon_k)=1$ or $d_k(\delta_k)=d_k(\delta_k+\varepsilon_k)=0$ then, to fulf\/il the
equation~\eqref{eqK}, one has to set $p_0(x)=1-d_k(x-k \theta)$ for $x\in[k\theta+\delta_k,k\theta+\delta_k+\varepsilon_k]$.
The function $p_k$ should then be def\/ined as previously by $\sqrt{d_k(x)(1-d_k(x))}$ for
$x\in[\delta_k,\delta_k+\varepsilon_k]$ and 0 elsewhere, so that~\eqref{eq0} remains fulf\/illed.
The only task to accomplish is to choose well the parameters~$\varepsilon_k$ and~$\delta_k$ to avoid the possible
intersection of~$d_k$ functions.
The parameters $\varepsilon_k$ should be such that the equations~\eqref{eps} remain satisf\/ied, and
$\delta_k=n\theta+\varepsilon_{k_1}+\cdots+\varepsilon_{k_s}$ for~$n$, $s \in \Z$ which depend on the concrete projection
one has chosen as a~starting point.
\end{proof}

\begin{figure}[t]\centering
\includegraphics{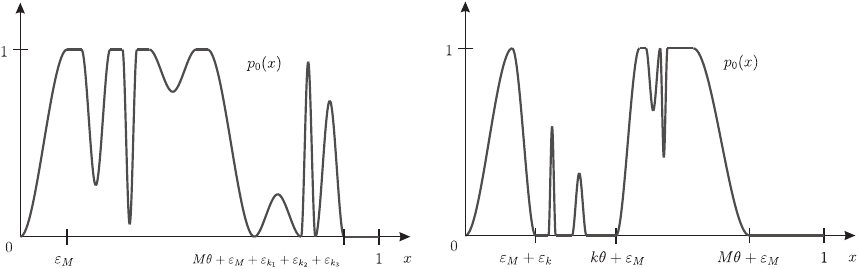}
\caption{Examples of $p_0$ functions for projections of traces $M\theta$ and $(M-k)\theta$.}\label{f4}
\end{figure}

Examples of $p_0$ functions of the described above projections are shown in Fig.~\ref{f4}.

By giving constructive proofs of Theorems~\ref{lm1} and~\ref{lm2} we have exhausted all of the possibilities of
constructing projections in $\NCA$ with the method described in Section~\ref{equations}.
To end this section let us note that the computation of the Chern number of the newly constructed projections does not
provide any new information.
Indeed, it is straightforward either from direct computations of the formula~\eqref{ChernSimp}, either from an
application of the results of~\cite{ConnesChern} that if we have a~projection~$p$ of trace $n \theta$, then $c_1(p)=n$.
In particular, the process of adding ``bump'' functions described in the proof of Theorem~\ref{lm2} does not change the
Chern class of a~projection.

\section{\texorpdfstring{Flip-symmetric projections}{Flip-symmetric projections}}
\label{symmetric}

The Powers--Rief\/fel type projection can be made invariant under the f\/lip automorphism
\begin{gather}
\label{flip}
\sigma \in \Aut A_{\theta}, \qquad \sigma(U) = U^{-1}, \qquad \sigma(V) = V^{-1}.
\end{gather}
This is accomplished by setting $\varepsilon_M = 1 - M \theta$ and requiring that $d_M(x) + d_M(\varepsilon_M - x) = 1$
for $x \in [0,\varepsilon_M]$ in formulae~\eqref{PRT0}--\eqref{PRTM}.
It is interesting to check for which of the more general projections presented in this paper the f\/lip symmetry can be
imposed\footnote{We are grateful to the anonymous referee for suggesting this interesting problem to us.}.

Requiring $\sigma(p) = p$ for projections of order $M$ translates to the following constraints on functions $p_k$:
\begin{gather}
\label{constraints}
p_0(x) = p_0(1-x),
\qquad
p_k(1-x) = \overline{p_k(x-k \theta)},
\qquad
\text{for}
\quad
k = 1, \ldots, M.
\end{gather}
The fact that $p_0$ is periodic with period 1 implies that for a~f\/lip-symmetric projection, $p_0$ should be symmetric
around $\frac{1}{2}$.

Let us f\/irst note that the f\/lip symmetry cannot be imposed on a~projection, which is ``glued'' from two or more segments
(see Fig.~\ref{f2}).
This is because the segments will necessarily have dif\/ferent lengths and thus $p_0$ cannot be symmetric around
$\frac{1}{2}$.
It is also clear (see Fig.~\ref{f4}) that the inclusion of a~``bump function'' breaks the f\/lip invariance of
a~projection.

\begin{figure}[t]\centering
\includegraphics{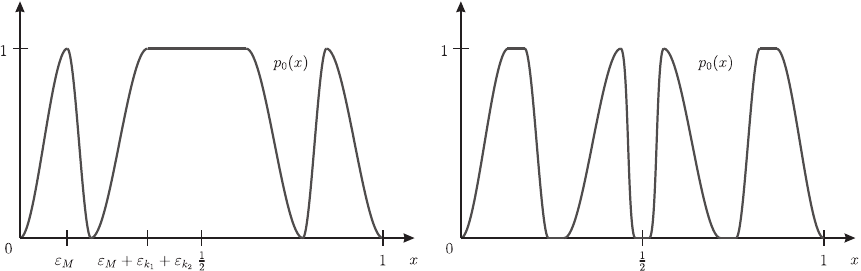}
\caption{Examples of $p_0$ functions for f\/lip-symmetric projections.}\label{f5}
\end{figure}

On the other hand, the ``cutting out'' procedure (see Fig.~\ref{f3}) can be performed in such a~way that the invariance
under the f\/lip automorphism is preserved.
Moreover, provided that the $\theta$ parameter is small enough, one can ``glue'' another Powers--Rief\/fel type projection
inside the ``cut-out'' region in a~symmetric way.
This process may be continued as long as the parameter~$\theta$ allows it.
However, this requires a~f\/ine tuning of~$\varepsilon_k$ parameters.
Examples of~$p_0$ functions of such f\/lip-symmetric projection are depicted in Fig.~\ref{f5}.

Let us now formulate the above considerations in a~precise way.
\begin{proposition}
\label{PropSym}
Let $k_j \in \{1, \ldots, M-1\}$ for $j = 1,2,\ldots,s$, $1 \leq s \leq M-1$ be such that
\begin{gather*}
M > k_1 > k_2 > \cdots > k_s > 0
\qquad
\text{and}
\qquad
\theta \in \left]c_1^-, c_1^+ \right] \cap \left[c_2^-, c_2^+ \right[
\cap \cdots \cap \left[c_s^-, c_s^+ \right] \setminus \{c_s^{\pm}\},
\end{gather*}
where $\pm = (-1)^s$ and
\begin{gather*}
c_j^+=
\begin{cases}\displaystyle
\frac{1}{2 (M-k_1+k_2-k_3+\cdots-k_j)}, & \text{for $j$ odd},
\vspace{1mm}\\
\displaystyle \frac{1}{2 (M-k_1+k_2-k_3+\cdots-k_{j-1})+k_j}, & \text{for $j$ even},
\end{cases}
\\
c_j^-=
\begin{cases}
\displaystyle \frac{1}{2 (M-k_1+k_2-k_3+\cdots+k_{j-1})-k_j}, & \text{for $j$ odd},
\vspace{1mm}\\
\displaystyle \frac{1}{2 (M-k_1+k_2-k_3+\cdots+k_j)}, & \text{for $j$ even}.
\end{cases}
\end{gather*}
Then, there exists a~projection $p \in \NCA$ of order $M$ invariant under the flip automorphism~\eqref{flip} $\sigma(p)=p$.
Moreover, $p$ is a~representant of the $K_0$ class $[(M-k_1+k_2-k_3+\cdots + (-1)^s k_s)\theta]$.
\end{proposition}
\begin{proof}
The construction of the projection $p$ goes as follows: We start with a~Powers--Rief\/fel type projection of order~$M$ \eqref{PRT0}--\eqref{PRTM}.
We make it f\/lip symmetric, i.e.~we set $\varepsilon_M = 1 - M \theta$ and require that $d_M$ be such that $d_M(x) +
d_M(\varepsilon_M - x) = 1$ for $x \in [0,\varepsilon_M]$.
Then, we ``cut out'' a~projection of trace $k_1 \theta$ (see formulae~\eqref{cut0}--\eqref{cutk}) also in
symmetry-preserving way.
This requires setting (compare the left plots in Figs.~\ref{f3} and~\ref{f5})
\begin{gather*}
\varepsilon_{k_1} = (2M-k_1) \theta - 1
\qquad
\text{and}
\qquad
d_{k_1}(x) + d_{k_1}(\varepsilon_{k_1} - x) = 1,
\qquad
\text{for}
\quad
x\in[\varepsilon_M,\varepsilon_M+\varepsilon_{k_1}].
\end{gather*}
This choice guarantees that constraints~\eqref{constraints} are fulf\/illed for both functions~$p_0$ and $p_{k_1}$.

Note however, that we need to have $0 < \varepsilon_{k_1} \leq k_1 \theta$ for $p$ to be a~projection, which is
equivalent to
\begin{gather*}
\frac{1}{2M-k_1} < \theta \leq \frac{1}{2(M-k_1)}
\quad
\Longleftrightarrow
\quad
\theta \in \left] c_1^-, c_1^+ \right].
\end{gather*}

Now, we ``glue'' a~Powers--Rief\/fel type projection of trace $k_2 \theta$ in the middle of the ``cut-out'' region.
To preserve the f\/lip symmetry we have to set
\begin{gather*}
\varepsilon_{k_2} = 1 - 2(M-k_1) \theta + k_2 \theta
\end{gather*}
and
\begin{gather*}
d_{k_2}(x) + d_{k_2}(\varepsilon_{k_2} - x) = 1,
\qquad
\text{for}
\quad
x\in[\varepsilon_M+\varepsilon_{k_1},\varepsilon_M+\varepsilon_{k_1}+\varepsilon_{k_2}].
\end{gather*}
But since $0 < \varepsilon_{k_2} \leq k_2 \theta$, the equation on $\varepsilon_{k_2}$ can be fulf\/illed only if
\begin{gather*}
\frac{1}{2(M-k_1+k_2)} \leq \theta < \frac{1}{2(M-k_1)+k_2}
\quad
\Longleftrightarrow
\quad
\theta \in \left[ c_2^-, c_2^+ \right[.
\end{gather*}

By performing further consecutive ``cut-outs'' and ``glueings'' we obtain the following conditions for all $j \in \{1,
\ldots, s\}$
\begin{gather*}
\varepsilon_{k_j} =
\begin{cases}
2 (M-k_1+k_2-k_3+\cdots-k_{j-1}) \theta - k_j \theta -1, & \text{for $j$ odd},
\\
1- 2 (M-k_1+k_2-k_3+\cdots-k_{j-1}) \theta - k_j, & \text{for $j$ even},
\end{cases}
\\
d_{k_j}(x) + d_{k_j}(\varepsilon_{k_j} - x) = 1,
\qquad
\text{for}
\quad
x \in [\varepsilon_M+\varepsilon_{k_1}+\cdots+\varepsilon_{k_{j-1}},\varepsilon_M+\varepsilon_{k_1}+\cdots +\varepsilon_{k_j}],
\end{gather*}
which can be met only if $\theta \in [c_j^{-},c_j^{+}] \setminus \{c_j^{\pm}\}$, with $\pm = (-1)^j$.

To conclude the proof let us remind the reader that the procedure of ``glueing'' a~projection of trace $k \theta$
increases the trace of the overall projection by $k\theta$ and ``cutting-out'' decreases it by~$k\theta$.
\end{proof}

We end this section with a~remark, that some of the presented symmetric projections cannot be constructed in any $\NCA$
with $\theta \in [0,1[$.
Let us take for instance $M=6$, $k_1 = 5$, $k_2 = 4$, $k_3 = 1$.
Then $c_1^- = \frac{1}{7}$, $c_1^+ = \frac{1}{2}$, $c_3^- = \frac{1}{9}$, $c_3^+ = \frac{1}{8}$, so
$]c_1^-,c_1^+]  \, \cap\,   ]c_3^-, c_3^+] = \varnothing$.

\section{Conclusion and open questions}
\label{conc}
Let us now summarise the obtained results and outline the directions of possible further investigations.

We have presented many projections, which generalise the standard Powers--Rief\/fel projection.
Some represented the same $K_0(\A_{\theta})$ class, but had dif\/ferent orders.
The others, conversely, had the same order, but dif\/ferent traces.
A natural question one can ask is what are the relations between the presented projections? The answer is provided by
Theorem~8.13 in~\cite{Rieffel3}.
It states that if two projections in $\NCA$ represent the same $K_0(\NCA)$ class (hence have the same trace), then, not
only they are unitarily equivalent in $M_{\infty}(\A_{\theta})$, but they are actually in the same path component of the
set of projections in $\NCA$ itself.
This means that there exists a~homotopy of projections in $\NCA$ for any two projections which have the same trace.
Indeed, if, for instance, one takes $d_{k}(t,x):= t d_{k}(x) + (1-t)$ instead of $d_k(x)$ with
$d_k(\delta_k)=d_k(\delta_k+\varepsilon_k)=1$ for the ``bump function'' used in the proof of Theorem~\ref{lm2}, then one
would obtain a~projection for all $t \in [0,1]$.

In consequence, from the topological point of view it is suf\/f\/icient to consider Powers--Rief\/fel type projections, since
they are the generators of the $K_0(\NCA)$ group.
On the other hand, the richness of the structure of projections may show up in applications.
In the proof of the Theorem~\ref{lm1} it has already been mentioned that the procedure of ``gluing'' the Powers--Rief\/fel
type projections is in fact equivalent to taking sums of mutually orthogonal projections.
However, the ``cutting out'' described subsequently does not admit an interpretation in terms of subtracting the
projections.
Indeed, it is straightforward to see that if one expresses a~projection $p$ of order~$M$ and trace $(M-k)\theta$ as $p=q-r$,
where~$q$ is a~Powers--Rief\/fel type projection of order~$M$, then $r$ would not be a~projection.
This shows that the newly found projections are not just linear combinations of the~$K_0$ generators.

The method adopted in this paper clearly does not pretend to cover every possible projection in~$\NCA$.
For a~most general projection in~$\NCA$ one would need to allow the order~$M$ of a~projection go to inf\/inity.
This would imply the need of working directly with the elements of the form~\eqref{a} and would require completely
dif\/ferent methods (see~\cite{Boca} for an example).

The puzzling thing about the newly found projections is that their existence in $\NCA$ seems to depend on the
noncommutativity parameter $\theta$ as stated in the theorems in Section~\ref{general}.
Unfortunately, the solutions presented there cannot be adapted to the case $n \theta > 1$, as it was done for the
Powers--Rief\/fel type projections in Section~\ref{sol}.
It is so because the translational symmetry of~\eqref{PRM1}--\eqref{PRM3}, used in the proof of Proposition~\ref{Thm} is
absent in general equations~\eqref{selfad}--\eqref{idemp2}.
Note that the discussed symmetry is also broken whenever we introduce the mentioned ``bump functions''.
What is more, the existence of f\/lip-symmetric projections in~$\NCA$ seems to depend on~$\theta$ in an even more peculiar
way.
Whether there is a~true dif\/ference in the structure of projections in~$\NCA$ depending on the noncommutativity parameter~$\theta$ or is it just an artefact of our method of solving the equations~\eqref{selfad}--\eqref{idemp2} remains an open question.

To conclude the paper, let us comment on the possible applications of the obtained results to the $D$-brane scenario in
Type II string theories.
As mentioned in the Introduction, projections in $\NCA$ correspond to solitonic f\/ield conf\/igurations which are
identif\/ied with $D$-branes~\cite{PRD,Krajewski,NFT}.
On one hand, unitarly equivalent projections yield gauge equivalent f\/ield conf\/igurations~\cite[Section~3.1]{NFT}, hence
the knowledge of $K_0(\NCA)$ alone seems to be suf\/f\/icient.
On the other hand, projections which cannot be written as linear combinations of $K_0$ generators provide
non-perturbative f\/ield conf\/igurations.
Moreover, the homotopy equivalence of projections may be exploited to study the soliton dynamics.
An example is provided in~\cite[Section~6.2]{NFT}, where the Boca projection~\cite{Boca}, which is homotopy equivalent
to the standard Powers--Rief\/fel projection, is used.
The possibility of adding ``bump functions'' to a~projection as described at the end of Section~\ref{general} indicates
the existence of an additional degree of freedom of the D-branes.
It would also be interesting to investigate the consequences for $D$-branes of a~projection being invariant under the
f\/lip symmetry.
Finally, let us note that the $D$-brane point of view suggests that the number of projections in $\NCA$ indeed depends
on the value of the deformation parameter~$\theta$ (see~\cite[Section~4]{Krajewski} or~\cite[Section~V]{PRD}).

\subsection*{Acknowledgements}

We would like to thank Andrzej Sitarz for his illuminating remarks.
Project operated within the Foundation for Polish Science IPP Programme ``Geometry and Topology in Physical Models''
co-f\/inanced by the EU European Regional Development Fund, Operational Program Innovative Economy 2007-2013.

\pdfbookmark[1]{References}{ref}
\LastPageEnding

\end{document}